\documentclass[twoside,12pt]{article}
\usepackage{amssymb,amsmath,bm,euscript}
\usepackage{graphicx,color,caption}
\usepackage{ifpdf}

\usepackage[colorlinks=true]{hyperref}

\usepackage{braket}

\usepackage{algorithm, algorithmic}

\newtheorem{proposition}{Proposition}[section]
\newtheorem{theorem}[proposition]{Theorem}

\newtheorem{corollary}[proposition]{Corollary}

\newcommand{\qed}{\hphantom{.}\hfill $\Box$\medbreak}

\def\S{\mathcal{S}}

\def\O{\mathcal{O}}
\def\R{\mathbb{R}}
\def\F{\mathcal{F}}
\def\I{\mathcal{I}}

\def\A{{\mathcal{A}}}
\def\B{\mathcal{B}}

\def\CC{\mathbb{C}}

\def\D{\mathcal{D}}
\def\U{\mathcal{U}}
\def\V{\mathcal{V}}
\def\S{\mathcal{S}}
\def\X{{\mathcal{X}}}

\def\Z{\mathcal{Z}}
\def\W{\mathcal{W}}
\def\O{{\mathcal{O}}}

\def\aa{{\bf a}}
\def\bb{{\bf b}}

\def\dd{{\bf d}}

\def\ee{{\bf e}}
\def\x{{\bf x}}

\def\s{{\bf s}}

\def\uu{{\bf u}}

\def\0{{\bf 0}}

\title{\bf{T-Quadratic Forms and Spectral Analysis of T-Symmetric Tensors}}
\author{ \hspace{1mm} Liqun Qi\thanks{
Department of Applied
    Mathematics, The Hong Kong Polytechnic University, Hung Hom,
    Kowloon, Hong Kong, China; ({\tt liqun.qi@polyu.edu.hk}).}
 \ and \
  Xinzhen Zhang\thanks{School of Mathematics, Tianjin University, Tianjin 300354 China; ({\tt xzzhang@tju.edu.cn}). This author's work was supported by NSFC (Grant No.  11871369). }
}

\begin{document}
\date{\today}
\maketitle

\begin{abstract}
An $n \times n \times p$ tensor is called a T-square tensor.   It arises from many applications, such as the image feature extraction problem and the multi-view clustering problem.  We may symmetrize a T-square tensor to a T-symmetric tensor.  For each T-square tensor, we define a T-quadratic form, whose variable is an $n \times p$ matrix, and whose value is a $p$-dimensional vector.   We define eigentuples and eigenmatrices for T-square tensors.   We show that a T-symmetric tensor has unique largest and smallest  eigentuples, and a T-quadratic form is positive semi-definite (definite) if and only if its smallest eigentuple is nonnegative (positive).  The relation between the eigen-decomposition of T-symmetric tensors, and the TSVD of general third order tensors are also studied.



\end{abstract}


\section{Introduction}

 We call a third order tensor $\A \in \Re^{n \times n \times p}$ a T-square tensor.  It was called an f-square tensor in \cite{MQW21}.   The representation tensor $\Z \in \Re^{n \times n \times p}$ arising in the multi-view clustering problem \cite{CXZ20} and the multi-view image feature extraction problem \cite{XCGZ21, XCGZ21a} is a T-square tensor.   Here $n$ is the number of the samples in the database, $p$ is the number of the views.

Suppose that $\A \in \Re^{n \times n \times p}$ is a T-square tensor.   Let $X \in \Re^{n \times p}$.   We may regard $X$ as a tensor $\X \in \Re^{n \times 1 \times p}$.   Define
\begin{equation} \label{e1.1}
F_\A(X) : = \X^\top * \A * \X,
\end{equation}
where $*$ is the T-product operation introduced in \cite{Br10, KBHH13, KM11}, and $\top$ is the transpose operation in the T-product sense.   In the next section, we will review the definition of T-product and its transpose concept.
We call $F_\A$ the T-quadratic form defined by $\A$.  Then for any $X \in \Re^{n \times p}$, $F_\A(X) \in \Re^p$.    If $F_\A(X) \ge \0$ for any $X \in \Re^{n \times p}$, then we say that the T-quadratic form $F_\A$ is T-positive semi-definite.   If $F_\A(X) > \0$ for any $X \in \Re^{n \times p}$, then we say that the T-quadratic form $F_\A$ is T-positive definite.

The T-positive semidefiniteness (definiteness) concept here is different from the T-positive semidefiniteness (definiteness) concept discussed in \cite{ZHW21}.  The T-positive semidefiniteness (definiteness) concept in \cite{ZHW21} is in the sense of nonnegative (positive) scalars.   Here, the concept is in the sense of  nonnegative (positive) vectors.   Thus, the T-positive semidefiniteness (definiteness) concept here is stronger and may reflect more correlative properties of T-square tensors.

The T-product operation, TSVD decomposition and tubal ranks were introduced by Kilmer and her collaborators in \cite{Br10, KBHH13, KM11}.  It is now widely used in engineering \cite{CXZ20, MQW20, MQW21, SHKM14, SNZ21, XCGZ21, XCGZ21a, YHHH16, ZHW21, ZSKA18, ZA17, ZEAHK14, ZLLZ18}.  In \cite{Br10}, Bradman defined real eigentuples and eigenmatrices for third order tensors in $\Re^{n \times n \times n}$.
Viewing the wide applications of T-product, TSVD decomposition and tubal ranks, the theory of eigentuples and eigenmatrices deserves to be further studied.

In this paper, we extend the concepts of eigentuples and eigenmatrices to T-square tensors and allow complex eigentuples and eigenmatrices.   We show that an $n \times n \times p$ T-symmetric tensor has unique largest eigentuple $\s_1 \in \Re^p$  and unique smallest eigentuple $\s_n \in \Re^p$
such that any real eigentuple $\s$ of $\A$ satisfies $\s_1 \ge \s \ge \s_n$.   We further show that a T-quadratic form is positive semidefinite (definite) if and only if the smallest eigentuple of the corresponding T-symmetric tensor is nonnegative (positive).

The T-quadratic function $F_\A$ maps $\Re^{n \times p}$ to $\Re^p$.   Its positive semidefiniteness (definiteness) requires $p$ quadratic polynomials of $np$ variables to be nonnegative (positive) simutaneously.  We
present its spectral conditions.   This theory is noval.

We then further study the relation between the eigen-decomposition of T-symmetric tensors, and the TSVD of general third order tensors.

The reset of this paper is distributed as follows.  We deliver some preliminary knowledge of T-product operations in the next section.   In Section 3, we define eigentuples and eigenmatrices for a T-square tensor, and show the existence of the largest and the smallest eigentuples of a T-symmetric tensor.   In Section 4, we prove that a T-symmetric tensor is positive semidefinite (definite) if and only if its smallest eigentuple is nonnegative (positive).   We study the relation between the eigen-decomposition of T-symmetric tensors, and the TSVD of general third order tensors in Section 5.   

\section{Preliminaries}

Let $\aa = (a_1, a_2, \cdots, a_p)^\top \in \CC^p$.   Then
$${\rm circ}(\aa):= \left(\begin{aligned} a_1\ & a_p & a_{p-1} & \cdots & a_2\ \\ a_2\ & a_1 & a_p\ & \cdots & a_3\ \\ \cdot\ \ \ & \ \cdot & \cdot\ \  & \cdots & \cdot\ \ \ \\
 \cdot\ \ \ & \ \cdot & \cdot\ \  & \cdots & \cdot\ \ \ \\
 a_p & a_{p-1} & a_{p-2} & \cdots & a_1 \end{aligned}\right),$$
 and circ$^{-1}($circ$(\aa)) := \aa$.

 Suppose that $\aa, \bb \in \CC^p$.   Define
 $$\aa \odot \bb = {\rm circ}(\aa)\bb.$$

 In \cite{KBHH13}, $\aa, \bb \in \Re^p$ are called tubal scalars.   Here, we extend them to $\CC^p$.

 In general, $\aa \odot \bb \not = \bb \odot \aa$.   We denote
 $$\aa^{\odot 2} := \aa \odot \aa.$$
 If $\aa \in \Re^p$ is nonnegative, then $\aa^{\odot 2}$ is also nonnegative.   However, if $\bb \in \Re^p$ is nonnegative, there may be no $\aa \in \Re^p$ such that $\aa^{\odot 2} = \bb$.   For example, let $p=2$,
 $\aa = (a_1, a_2)^\top$, $\bb = (b_1, b_2)^\top$ and $\aa^{\odot 2} = \bb$.   Then we have
 $b_1 = a_1^2 + a_2^2$ and $b_2 = 2a_1a_2$.    To satisfy these two equations, we must have $b_1 \ge b_2$.
 We say that $\bb \in \Re^p$ is a square tubal scalar if it is nonnegative and there is an $\aa \in \Re^p$, such that $\aa$ is nonnegative and $\aa^{\odot 2} = \bb$.

 For $\aa = (a_1, \cdots, a_p)^\top \in \Re^p$, denote $|\aa| := (|a_1|, \cdots, |a_p|)^\top$.

 {\bf Question 1} Suppose that $\bb \in \Re^p$ is a square tubal scalar.   Is there a unique $\aa \in \Re^p$, such that $\aa$ is nonnegative and $\aa^{\odot 2} = \bb$?

 \begin{proposition} \label{p2.1}
 $(\CC^p, +, \odot)$ is a commutative ring with unity $\ee = (1, 0, \cdots, 0)^\top \in \CC^p$, where $+$ is the vector addition.
 \end{proposition}

 Proposition \ref{p2.1} extends Theorem 3.2 of \cite{Br10} from $\Re^p$ to $\CC^p$, as we need to consider
 complex eigentuples for third order real tensors.  The proof is almost the same.   Hence, we omit the proof.

 Note that the operation $\odot$ is different from vector convolution.  For  $\aa, \bb \in \CC^p$, the vector convolution of $\aa$ and $\bb$ is in $\CC^{2p-1}$.

For $X \in \CC^{n \times p}$ and $\aa \in \CC^p$, define
$$\aa \circ X = X{\rm circ}(\aa).$$

\begin{proposition} \label{p2.2}
Let $\aa, \bb \in \CC^p$, and $X, Y \in \CC^{n \times p}$.   Then

1. $\aa \circ (X+Y) = \aa \circ X + \aa \circ Y$;

2. $(\aa + \bb) \circ X = \aa \circ X + \bb \circ X$;

3. $\aa \circ (\bb \circ X) = (\aa \odot \bb) \circ X$;

4. Let $\ee = (1, 0, \cdots, 0)^\top \in \CC^p$ as in Proposition \ref{p2.1}.  Then $\ee \circ X = X$ for all $X \in \CC^{n \times p}$.   Furthermore, $\ee$ is the unique element in $\CC^p$ with this property.
\end{proposition}
{\bf Proof}  This proposition extends Theorem 3.5 of \cite{Br10} from $\CC^{p \times p}$ to $\CC^{n \times p}$, except the second half of item 4 is additional.   The proof of the other part except the second half of item 4 is almost the same as the proof of Theorem 3.5 of \cite{Br10}.  We omit this part and now prove the second half of item 4.  Suppose $\aa \circ X = X$ for all $X \in \CC^{n \times p}$.  Then
$X {\rm circ}(\aa) = X$ for all $X \in \CC^{n \times p}$.  This implies ${\rm circ}(\aa) = I_p$, the identity matrix of $\Re^{p \times p}$.   Thus, $\aa = {\rm circ}^{-1}(I_p) = \ee$.
\qed

For a third order tensor $\A \in \Re^{m \times n \times p}$, its frontal slices are denoted as $A^{(1)}, \cdots, A^{(p)} \in \Re^{m \times n}$.   As in  \cite{Br10,KBHH13,KM11}, define
$${\rm bcirc}(\A):= \left(\begin{aligned} A^{(1)}\ & A^{(p)} & A^{(p-1)} & \cdots & A^{(2)}\ \\ A^{(2)} & A^{(1)} & A^{(p)} & \cdots & A^{(3)}\\ \cdot\ \ \ & \ \cdot & \cdot\ \  & \cdots & \cdot\ \ \ \\
 \cdot\ \ \ & \ \cdot & \cdot\ \  & \cdots & \cdot\ \ \ \\
 A^{(p)} & A^{(p-1)} & A^{(p-2)} & \cdots & A^{(1)} \end{aligned}\right),$$
 and bcirc$^{-1}($bcirc$(\A)):= \A$.

Various T-product structured properties of third order tensors are based upon their block circulant matrix versions.   For a third order tensor $\A \in \Re^{m \times n \times p}$, its transpose
 can be defined as
 $$\A^\top = {\rm bcirc}^{-1}[({\rm birc}(\A))^\top].$$
 This will be the same as the definition in \cite{Br10, KBHH13, KM11}.  The identity tensor $\I_{nnp}$ may also be defined as
 $$\I_{nnp} = {\rm bcirc}^{-1}(I_{np}),$$
 where $I_{np}$ is the identity matrix in $\Re^{np \times np}$.  

 However, a third order tensor $\S$ in $\Re^{m \times n \times p}$ is f-diagonal in the sense of \cite{Br10, KBHH13, KM11} if all of its frontal slices $S^{(1)}, \cdots, S^{(p)}$ are diagonal.   In this case, bcirc$(\S)$ may not be diagonal.

  For a third order tensor $\A \in \Re^{m \times n \times p}$, it is defined \cite{Br10, KM11} that
 $${\rm unfold}(\A) := \left(\begin{aligned} A^{(1)}\\ A^{(2)}\\ \cdot\ \ \\ \cdot\ \ \\ \cdot\ \ \\ A^{(p)}\end{aligned}\right) \in \Re^{mp \times n},$$
and fold$($unfold$(\A)) := \A$.   For $\A \in \Re^{m \times s \times p}$ and $\B \in \Re^{s \times n \times p}$, the T-product of $\A$ and $\B$ is defined as
$\A * \B :=$ fold$(${bcirc$(\A)$unfold$(\B) \in \Re^{m \times n \times p}$.   Then, we see that
$$\A * \B = {\rm bcirc}^{-1}({\rm bcirc}(\A){\rm bcirc}(\B)).$$
Thus, the bcirc and bcirc$^{-1}$ operations not only form a one-to-one relationship between third order tensors and block circulant matrices, but their product operation is reserved.

 \underline{The Standard Form of a Real f-Diagonal Tensor.}   Let $\S = (s_{ijk}) \in \Re^{m \times n \times p}$  be a f-diaginal tensor.   Let $\s_j = (s_{jj1}, \s_{jj2}, \cdots, s_{jjp})^\top$ be the $jj$th tube of $\S$ for $j = 1, \cdots, \min \{m, n \}$.   We say that $\S$ is in its standard form if $\s_1 \ge \s_2 \ge \cdots \ge \s_{\min \{m, n\}}$.

\section{Eigentuples and Eigenmatrices of T-Square Tensors}


For a matrix $X \in \CC^{n \times p}$, let its column vectors be $\x^{(1)}, \cdots, \x^{(p)}$.
Define
$${\rm unfold}(X) := \left(\begin{aligned} \x^{(1)}\\ \x^{(2)}\\ \cdot\ \ \\ \cdot\ \ \\ \cdot\ \ \\ \x^{(p)}\end{aligned}\right) \in \CC^{np},$$
and fold$($unfold$(X)) : = X$.       Then we define the T-product of $\A$ and $X$ as
 $$\A * X = {\rm fold}({\rm bcirc}(\A){\rm unfold}(X)).$$

 Thus, $\A * X \in \CC^{m \times p}$.

 \bigskip

 \bigskip

 We now define eigentuples and eigenmatrices of T-square tensors.   Suppose that $\A \in \Re^{n \times n \times p}$ is a T-square tensor, $X \in \CC^{n \times p}$ and $X \not = O$, $\dd \in \CC^p$, and
 \begin{equation} \label{e1}
 \A * X = \dd \circ X.
 \end{equation}
 Then we call $\dd$ an eigentuple of $\A$, and $X$ an eigenmatrix of $\A$, corresponding to the eigentuple $\dd$.

 The eigentuple and eigenmatrix concepts extend the eigentuple and eigenmatrix concepts of \cite{Br10} from $\Re^{p \times p \times p}$ to $\Re^{n \times n \times p}$ and allow complex eigentuples and  eigenmatrices.


We aim to study T-positive semi-definiteness and T-positive definiteness of the T-quadratic form $F_\A$, defined in (\ref{e1.1}).    This would not be easy by using the eigentuples of $\A$, as even for real square matrices,  their eigenvalues may not be real.   Thus, as in the matrix case, we symmetrize the T-square tensor $\A$.

Let $\A \in \Re^{n \times n \times p}$ be a T-square tensor.    We say that $\A$ is T-symmetric if
$\A = \A^\top$.  T-symmetric tensors have been studied in \cite{ZHW21}.   We have the following propositions.

\begin{proposition}
Suppose that $\A \in \Re^{n \times n \times p}$.   Then $\A + \A^\top$ is a T-symmetric tensor.  Then
$\A$ is positive semidefinite (definite) if and only if the T-symmetric tensor $\A + \A^\top$ is positive semidefinite (definite).
\end{proposition}
{\bf Proof} Since
$$\left(\A + \A^\top\right)^\top = \A^\top + \left(\A^\top\right)^\top = \A + \A^\top,$$
$\A + \A^\top$ is T-symmetric.

For $X \in \Re^{n \times p}$, regard it as a tensor $\X \in \Re^{n \times 1 \times p}$.  We have
$$F_\A(X) = \X^\top * \A * \X = (\X^\top * \A * \X)^\top = \X^\top * \A^\top * \X = {1 \over 2}
\X^\top * (\A + \A^\top) * \X.$$
Thus, $\A$ is positive semidefinite (definite) if and only if the T-symmetric tensor $\A + \A^\top$ is positive semidefinite (definite).
\qed

We thus study the eigentuples of T-symmetric tensors, and use them to analyze positive semidefiniteness (definiteness) of these tensors.

The following proposition holds obviously.

\begin{proposition}
A T-square tensor $\A \in \Re^{n \times n \times p}$ is T-symmetric if and only if bcirc$(\A)$ is symmetric.  A T-square tensor $\A \in \Re^{n \times n \times p}$ is invertible if and only if bcirc$(A)$ is invertible.   In this case, we have
$$\A^{-1} = {\rm bcirc}^{-1}({\rm bcirc}(\A^{-1})).$$
Furthermore, $\A$ is orthogonal in the sense of \cite{Br10, KBHH13, KM11} if and only if bcirc$(\A)$ is orthogonal.
\end{proposition}

We have the following theorem.

\begin{theorem} \label{t3.2}
Suppose that $\A \in \Re^{n \times n \times p}$ is a T-symmetric tensor.   Then there are orthogonal tensor $\U \in \Re^{n \times n \times p}$ and T-symmetric f-diagonal tensor $\D \in \Re^{n \times n \times p}$ such that
\begin{equation} \label{e3.2}
\A = \U * \D * \U^\top.
\end{equation}
Let the frontal slices of $\D$ be $D^{(1)}, \cdots, D^{(p)}$.  If $\hat \D \in \Re^{n \times n \times p}$
is another T-symmetric f-diagonal tensor, whose frontal slices $\hat D^{(1)}, \cdots, \hat D^{(p)}$ are resulted from switching some diagonal elements of $D^{(1)}, \cdots, D^{(p)}$, then there is an  orthogonal tensor $\hat \U \in \Re^{n \times n \times p}$, such that
\begin{equation} \label{e3.3}
\A = \hat \U * \hat \D * \hat \U^\top.
\end{equation}

\end{theorem}
{\bf Proof}   Block circulant matrices can be block diagonalized with normalized discrete Fourier transformation (DFT) matrix, which is unitary.   Then, as in (3.1) of \cite{KM11}, we have
\begin{equation} \label{e3.4}
(F_p \otimes I_n) \cdot {\rm bcirc}(\A) \cdot (F_p^* \otimes I_n) = {\rm diag}(D_1, \cdots, D_p),
\end{equation}
where $F_p$ is the $p \times p$ DFT matrix, $F_p^*$ is its conjugate transpose, $\cdot$ is the standard matrix multiplication, $\otimes$ denotes the Kronecker product.    Since bcirc$(\A)$ is symmetric, by taking conjugate transpose of (\ref{e3.4}), we see that $D_1, \cdots, D_p$ in (3.1) of \cite{KM11} are all hermite.   Applying the eigen-decomposition of $D_i = U_i\Sigma_iU_i^\top$ for $i = 1, \cdots, p$, we have
\begin{equation} \label{e3.5}
{\rm diag}(D_1, \cdots, D_p) = {\rm diag}(U_1, \cdots, U_p){\rm diag}(\Sigma_1, \cdots, \Sigma_p){\rm diag}(U_1^\top, \cdots, U_p^\top).
\end{equation}
Apply $(F_p^* \otimes I_n)$ to the left and $(F_p \otimes I_n)$ to the right of each of the block diagonal matrices in (\ref{e3.5}).  In each of the three cases, the resulting triple product results in a block circulant matrix.   We have
$${\rm bcirc}(\A) = {\rm bcirc}(\U) {\rm bcirc}(\D){\rm bcirc}(\U^\top).$$
This implies (\ref{e3.2}).  Then we have
$$\D = \U^\top * \A * \U,$$
and
$$\D^\top = \left(\U^\top * \A * \U\right)^\top = \U^\top * \A^\top * \U = \U^\top * \A * \U = \D,$$
as $\A$ is T-symmetric.  Thus, $\D$ is also T-symmetric.

Switching the order of eigenvalues in the eigen-decomposition $D_i = U_i\Sigma_iU_i^\top$ for $i = 1, \cdots, p$, we have (\ref{e3.3}).
\qed

We call (\ref{e3.2}) a T-eigen-decomposition (TED) of $\A$.

\begin{corollary}
Suppose that $\A \in \Re^{n \times n \times p}$ is a T-symmetric tensor and (\ref{e3.2}) holds.   Denote $\A^{*2} = \A * \A$ and $\A^{*k} = \A^{*(k-1)} * \A$ for any integer $k \ge 3$.  Then for any positive integer $k$, $\A^{*k}$ is still T-symmetric, and we have
$$\A^{*k} =  \U * \D^{*k} * \U^\top.$$
\end{corollary}

\begin{corollary}
Suppose that $\A \in \Re^{n \times n \times p}$ is a T-symmetric tensor and (\ref{e3.2}) holds.   Then $\A^{-1}$ exists if and only if $\D^{-1}$ exists.   If they exist, then they are T-symmetric and
$$\A^{-1} =  \U *  \D^{-1} * \U^\top.$$
\end{corollary}

We may rewrite (\ref{e3.2}) as
\begin{equation} \label{e3.6}
\A * \U = \U * \D,
\end{equation}
or
\begin{equation} \label{e3.7}
{\rm bcirc}(\A){\rm bcirc}(\U) = {\rm bcirc}(\U){\rm bcirc}(\D),
\end{equation}
Denote the $j$th lateral slice of $\U$ by $U_j \in \Re^{n \times p}$ for $j = 1, \cdots, n$.
Consider the $j$th column of (\ref{e3.7}) for $j = 1, \cdots, n$.   Let $\D = (d_{ijk})$.  Then
$d_{ijk}=0$ if $i \not = j$.  Let $d_{11k} \ge d_{22k} \ge \cdots \ge d_{nnk}$ for $k =1, \cdots, p$.       We have
\begin{equation} \label{e3.8}
\A * U_j = \dd_j \circ U_j,
\end{equation}
where $\dd_j = (d_{jj1}, d_{jjp}, d_{jj(p-1)}, \cdots, d_{jj2})^\top$.  Since $\U$ is orthogonal, $U_j \not = O$.  Thus, $\dd_j$ is an eigentuple of $\A$ with an eigenmatrix $U_j$.

For a matrix $U \in \Re^{n \times p}$, let its column vectors are $\uu^{(1)}, \cdots, \uu^{(p)}$.  Then
$$U = (\uu^{(1)}, \uu^{(2)}, \cdots, \uu^{(p)}).$$
Denote $U^{[0]} = U$,
$$U^{[1]} = (\uu^{(p)}, \uu^{(1)}, \cdots, \uu^{(p-1)}),$$
$$U^{[2]} = (\uu^{(p-1)}, \uu^{(2)}, \cdots, \uu^{(p-1)}),$$
$$\cdots,$$
$$U^{[p-1]} = (\uu^{(2)}, \uu^{(3)}, \cdots, \uu^{(1)}).$$
Consider the $(n+j)$th column of (\ref{e3.7}).   We have
\begin{equation} \label{e3.8}
\A * U_j^{[1]} = \dd_j \circ U_j^{[1]}.
\end{equation}
Thus, $U_j^{[1]}$ is also an eigenmatrix of $\A$, associated with the eigentuple $\dd_j$.
Similarly, $U_j^{[2]}, \cdots, U_j^{[p-1]}$ are also eigenmatrices of $\A$, associated with the eigentuple $\dd_j$.

Consider the set of eigenmatrices
$$T = \left\{ U_j^{[k]} : j =1, \cdots, n, k = 0, \cdots, p-1 \right\}.$$
Then $T$ forms an orthonormal basis of $\Re^{n \times p}$.  For any two distinct members $W, V \in T$,
let $\W$ and $\V$ be the corresponding $n \times 1 \times p$ tensors.  Then we have
\begin{equation} \label{e3.11}
\W^\top * \W = \I_{11p},
\end{equation}
and
\begin{equation} \label{e3.12}
\W^\top * \V = \O_{11p}.
\end{equation}

Viewing (\ref{e3.3}), we may switch the order in $\{ d_{11k}, d_{22k}, \cdots, d_{nnk} \}$ for any $k =1, \cdots, p$.  The resulted $\hat \dd_j$, $j = 1, \cdots, n$ are still eigentuples of $\A$.   Hence, the number of eigentuples of $\A$ is large.   But we may always take $\dd_1, \cdots, \dd_n$ in its standard form.

Combining the orthogonality of $\U$, we have the following theorem.

\begin{theorem} \label{t3.4}
Suppose that $\A \in \Re^{n \times n \times p}$ is a T-symmetric tensor.   Then $\A$ has real eigentuples
$\dd_1, \cdots, \dd_n$, such that $\dd_1 \ge \dd_2 \cdots \ge \dd_n$.   For each $j$, $j = 1, \cdots, n$, there are real eigenmatrices
$U_j^{[0]}, \cdots, U_j^{[p-1]}$, of $\A$, associated with the eigentuple $\dd_j$.   These
$np$ eigenmatrices form an orthonormal basis of $\Re^{n \times p}$.

\end{theorem}

We call the eigentuples $\{ \dd_1, \cdots, \dd_n \}$, satisfying $\dd_1 \ge \dd_2 \cdots \ge \dd_n$, in Theorem \ref{t3.4} the set of the principal eigentuples of $\A$. 

If $\A = \I_{nnp}$, then $\U = \D = \I_{nnp}$.   Therefore, $\dd_1 = \cdots = \dd_n = (1, 0, \cdots, 0)^\top$.  If $\A$ has a set of principal eigentuples $\dd_j = (d_{j1}, \cdots, d_{jp})^\top$ for $j=1, \cdots, n$, then $\A + \lambda \I_{nnp}$ has a set of principal eigentuples $\dd_j = (d_{j1}+\lambda, \cdots, d_{jp}+\lambda)^\top$ for $j=1, \cdots, n$.


We are not sure whether all eigentuples of a T-symmetric tensor are real and two eigenmatrices associated with two distinct eigentuples of a T-symmetric tensor are orthogonal to each other.  However,
but we can prove the following theorem.

\begin{theorem} \label{t3.5}
Suppose that $\A \in \Re^{n \times n \times p}$ is a T-symmetric tensor, and $\{ \dd_1, \cdots, \dd_n \}$
is a set of principal eigentuples of $\A$ such that $\dd_1 \ge \cdots \ge \dd_n$.
Then for any real eigentuple $\dd_0$ of $\A$, we have
\begin{equation} \label{e3.13}
\dd_1 \ge \dd_0 \ge \dd_n.
\end{equation}
\end{theorem}
{\bf Proof}   Assume that there is an eigenmatrix $V \in \CC^{n \times p}$ such that
$$\A * V = \dd_0 \circ V.$$
Taking conjugate, we have
$$\A * \bar V = \dd_0 \circ \bar V.$$
Let $W = V + \bar V$.  Then $W$ is real and
$$\A * W = \dd_0 \circ W.$$
If $W$ is nonzero, then $U$ is a real eigenmatrix of $\A$, associated with $\dd_0$.   Otherwise, $V$ is pure imaginary.   Letting $\hat W = \sqrt{-1} V$, we still have a real eigenmatrix of $\A$, associated with $\dd_0$.    Without loss of generality, assume $W$ is such a real eigenmatrix.

Let $U_j^{[0]}, \cdots, U_j^{[p-1]}$ be the eigenmatrices of $\A$ in Theorem \ref{t3.4}.   Then we have
real coefficients $\alpha_j^{[0]}, \cdots, \alpha_j^{[p-1]}$, for $j = 1, \cdots, n$, such that
$$W = \sum_{j=1}^n \sum_{k=1}^p \alpha_j^{[k-1]}U_j^{[k-1]}.$$
Let $\U_j^{[k]}$ be the $n \times 1 \times p$ tensors corresponding to $U_j^{[k]}$ for $j = 1, \cdots, n$ and $k = 0, \cdots, p-1$.   Let $\W$ be the $n \times 1 \times p$ tensor corresponding to $W$.   Let $\D_j$ be the $1 \times 1 \times p$ tensors corresponding to $\dd_j$ for $j = 0, \cdots, n$.
Then
\begin{eqnarray*}
& & \W^\top * \A * \W \\
& = & \W^\top * \W * \D_0 \\
& = & \left(\sum_{j=1}^n \sum_{k=1}^p \alpha_j^{[k-1]} \U_j^{[k-1]}\right)^\top * \left(\sum_{j=1}^n \sum_{k=1}^p \alpha_j^{[k-1]} \U_j^{[k-1]}\right) * \D_0 \\
& = & \sum_{j=1}^n \left(\sum_{k=1}^p \alpha_j^{[k-1]}\right)^2 \I_{11p} * \D_0 \\
& = & \sum_{j=1}^n \left(\sum_{k=1}^p \alpha_j^{[k-1]}\right)^2 \D_0.
\end{eqnarray*}
On the other hand,
\begin{eqnarray*}
& & \W^\top * \A * \W \\
& = & \W^\top * \A * \left(\sum_{j=1}^n \sum_{k=1}^p \alpha_j^{[k-1]} \U_j^{[k-1]}\right) \\
& = & \left(\sum_{j=1}^n \sum_{k=1}^p \alpha_j^{[k-1]} \U_j^{[k-1]}\right)^\top * \left(\sum_{j=1}^n \sum_{k=1}^p \alpha_j^{[k-1]} \A * \U_j^{[k-1]}\right)  \\
& = & \left(\sum_{j=1}^n \sum_{k=1}^p \alpha_j^{[k-1]} \U_j^{[k-1]}\right)^\top * \left(\sum_{j=1}^n \sum_{k=1}^p \alpha_j^{[k-1]} \U_j^{[k-1]} * \D_j \right)  \\
& = & \sum_{j=1}^n \left(\sum_{k=1}^p \alpha_j^{[k-1]}\right)^2 \I_{11p} * \D_j \\
& = & \sum_{j=1}^n \left(\sum_{k=1}^p \alpha_j^{[k-1]}\right)^2 \D_j.
\end{eqnarray*}

From this, we have
$$\sum_{j=1}^n \left(\sum_{k=1}^p \alpha_j^{[k-1]}\right)^2 \D_0 = \sum_{j=1}^n \left(\sum_{k=1}^p \alpha_j^{[k-1]}\right)^2 \D_j,$$
i.e.,
$$\sum_{j=1}^n \left(\sum_{k=1}^p \alpha_j^{[k-1]}\right)^2 \dd_0 = \sum_{j=1}^n \left(\sum_{k=1}^p \alpha_j^{[k-1]}\right)^2 \dd_j.$$
The inequality (\ref{e3.13}) is obtained.
\qed

\begin{corollary} \label{c3.6}
The eigentuples $\dd_1$ and $\dd_n$ are unique to $\A$.
\end{corollary}

We call $\dd_1$ the largest eigentuple of $\A$, and $\dd_n$ the smallest eigentuple of $\A$.


\section{T-Symmetric Positive Semidefinite Tensors}

Suppose that $\A \in \Re^{n \times n \times p}$ is a T-square tensor.   Then by Proposition \ref{p3.1},
the T-quadratic form is positive semidefinite (definite) if and only if the T-symmetric tensor $\A + \A^\top$ is positive semidefinite (definite).   This stimulates us to study positive semidefiniteness (definiteness) of a T-symmetric tensor.


\begin{theorem} \label{t4.1}
Suppose that $\A \in \Re^{n \times n \times p}$ is a T-symmetric tensor and it has a set of principal eigentuples
$\dd_1, \cdots, \dd_n$, such that $\dd_1 \ge \dd_2 \cdots \ge \dd_n$.
Then $\A$ is positive semidefinite (definite) if and only if the smallest eigentuple $\dd_n \ge (>) \0$.
\end{theorem}
{\bf Proof} By Theorem \ref{t3.4}, $\dd_1 \ge \dd_2 \cdots \ge \dd_n$, and for each $j$, $j = 1, \cdots, n$, there are real eigenmatrices
$U_j^{[0]}, \cdots, U_j^{[p-1]}$, of $\A$, associated with the eigentuple $\dd_j$, such that these
$np$ eigenmatrices form an orthonormal basis of $\Re^{n \times p}$.

If $\dd_n$ is not nonnegative, let $U = U_n^{[0]}$ and $\U_n^{[0]}$ be the corresponding $n \times 1 \times p$ tensor.   Let $\L$ be the $1 \times 1 \times p$ tensor corresponding to $\dd_n$.    Then
$$F_\A(U) = \left(\U_n^{[0]}\right)^\top * \A * \U_n^{[0]} =  \left(\U_n^{[0]}\right)^\top * \U_n^{[0]} * \L = \I_{nnp} * \L \not \ge \0,$$
which implies that $F_\A$ is not positive semi-definite.   Similarly, if $\dd_n$ is not positive, then $F_\A$ is not positive definite.

On the other hand, suppose that $\dd_n \ge 0$.  Let $\U_j^{[k]}$ be the  $n \times 1 \times p$ tensor corresponding to $U_j^{[k]}$ for $j=1, \cdots, n$ and $k = 0, \cdots, p-1$.  Let $X \in \Re^{n \times p}$.
Then there are real coefficients $\alpha_j^{[k]}$ for $j=1, \cdots, n$ and $k = 0, \cdots, p-1$, such that
$$X = \sum_{j=1}^n \sum_{k=0}^{p-1} \alpha_j^{[k]}U_j^{[k]}.$$
Let $\L_j$ be the $1 \times 1 \times p$ tensor corresponding to $\s_j$ for $j = 1, \cdots, n$.
We have
\begin{eqnarray*}
F_\A(X) & = & \left(\sum_{i=1}^n \sum_{l=0}^{p-1} \alpha_i^{[l]}\U_i^{[l]}\right)^\top * \A * \left(\sum_{j=1}^n \sum_{k=0}^{p-1} \alpha_j^{[k]}\U_j^{[k]}\right)\\
& = & \sum_{i=1}^n \sum_{j=1}^n \sum_{l=0}^{p-1} \sum_{k=0}^{p-1} \alpha_i^{[l]}\alpha_j^{[k]}
(\U_i^{[l]})^\top * \A * \U_j^{[k]}\\
& = & \sum_{i=1}^n \sum_{j=1}^n \sum_{l=0}^{p-1} \sum_{k=0}^{p-1} \alpha_i^{[l]}\alpha_j^{[k]}
(\U_i^{[l]})^\top * \U_j^{[k]} * \L_j\\
& = & \sum_{j=1}^n \sum_{k=0}^{p-1} \left(\alpha_j^{[k]}\right)^2\dd_j\\
& \ge & \0.
\end{eqnarray*}
Thus, $F_\A$ is positive semidefinite.   Similarly, if $\dd_n \ge \0$, then $\F_\A$ is positive definite.
\qed

\section{Relation with TSVD of General Third Order Tensors}

Suppose that $\A \in \Re^{m \times n \times p}$.  By \cite{KM11}, $\A$ has a T-singular value decomposition (TSVD):
\begin{equation} \label{e5.14}
\A = \U * \S * \V^\top,
\end{equation}
where $\U \in \Re^{m \times m \times p}$ and $\V \in \Re^{n \times n \times p}$ are orthogonal tensors, $\S \in \Re^{m \times n \times p}$ is a f-diagonal tensor.

\begin{theorem} \label{t5.1}
Suppose that $\A \in \Re^{m \times n \times p}$ with TSVD (\ref{e5.14}).   Then $\A * \A^\top \in \Re^{m \times m \times p}$ and $\A^\top * \A \in \Re^{n \times n \times p}$ are T-symmetric positive semi-definite tensors with TED
$$\A * \A^\top = \U * (\S * \S^\top) * \U^\top,$$
and
$$\A^\top * \A = \V * (\S^\top * \S) * \V^\top,$$
respectively.
\end{theorem}


 We now define singular tuples and singular matrices of general third order tensors.   Suppose that $\A \in \Re^{m \times n \times p}$ is a third order tensor, $X \in \Re^{n \times p}$, $X \not = O$, $Y \in \Re^{m \times p}$, $Y \not = O$, $\s \in \Re^p$, and
 \begin{equation} \label{e5.15}
 \A * X = \s \circ Y
 \end{equation}
 and
  \begin{equation} \label{e5.16}
 \A^\top * Y = \s \circ X.
 \end{equation}
 Then we call $\s$ a singular tuple of $\A$, $X$ a right singular matrix of $\A$,  $Y$ a right singular matrix of $\A$, corresponding to the singular tuple $\s$.

 \begin{theorem} \label{t5.2}
Suppose that $\A \in \Re^{m \times n \times p}$ is a third order tensor.   Without loss of generality, assume that $n \le m$.  Then $\A$ has singular tuples
$\s_1 \ge \s_2 \ge \cdots \ge \s_n \ge \0$.   For each $j$, $j = 1, \cdots, n$, there are right singular matrices
$U_j^{[0]}, \cdots, U_j^{[p-1]}$, and left singular matrices
$V_j^{[0]}, \cdots, V_j^{[p-1]}$, of $\A$, associated with the singular tuple $\s_j$.   The $np$ singular matrices $U_j^{[0]}, \cdots, U_j^{[p-1]}$ form an orthonormal basis of $\Re^{n \times p}$, and the $np$ singular matrices $V_j^{[0]}, \cdots, V_j^{[p-1]}$ form a part of an orthonormal basis of $\Re^{m \times p}$, respectively.

Furthermore, $\A^\top * \A \in \Re^{n \times n \times p}$ and $\A * \A^\top \in \Re^{m \times m \times p}$
are two T-symmetric positive semidefinite tensors.   The tensor $\A^\top * \A$ has $n$ nonnegative eigentuples $\s_1^{\odot 2}, \cdots, \s_n^{\odot 2}$.  For each $j$, $j = 1, \cdots, n$, there are real eigenmatrices $U_j^{[0]}, \cdots, U_j^{[p-1]}$, of $\A$, associated with the eigentuple $\s_j^{\odot 2}$.
The tensor $\A * \A^\top$ has $n$ nonnegative eigentuples $\s_1^{\odot 2}, \cdots, \s_n^{\odot 2}$.  For each $j$, $j = 1, \cdots, n$, there are real eigenmatrices $V_j^{[0]}, \cdots, V_j^{[p-1]}$, of $\A$, associated with the eigentuples $\s_j^{\odot 2}$.  If $n < m$, for each $j$, $j = n+1, \cdots, m$, there are real eigenmatrices $V_j^{[0]}, \cdots, V_j^{[p-1]}$, of $\A$, associated with the zero eigentuple $\0 \in \R^p$.  The $mp$ singular matrices $V_j^{[0]}, \cdots, V_j^{[p-1]}$ form an orthonormal basis of $\Re^{m \times p}$.
\end{theorem}

\bigskip

{\bf Acknowledgment}  We are thankful to Prof. Yicong Zhou and Dr. Dongdong Chen for the discussion on the multi-view clustering problem and the image feature extraction problem.

\end{document}